\documentclass[a4paper,12pt]{article}
\begin{document}

\title{Differentiations of operator algebras over non-archimedean fields.}
\author{Ludkovsky S.V.}
\date{15 September 2012}
\maketitle
\begin{abstract}
 Differentiations of operator algebras over non-archimedean
spherically complete fields are investigated. Theorems about a
differentiation being internal are demonstrated.\footnote{ 2010
Mathematics subject classification: 47A10, 47A25, 47L10
\par keywords: operator, algebra, field, differentiation }
\end{abstract}

\par address: Department of Applied Mathematics,
\par Moscow State Technical University MIREA,
av. Vernadsky 78,
\par Moscow, Russia
\par e-mail: Ludkowski@mirea.ru

\section{Introduction.}
Differentiations of operator algebras over the complex field were
investigated in \cite{kadannm66,sakaiannm66,losannm2008}. It was
shown that derivations of $C^*$-algebras and von Neumann algebras
are internal. But the case of operator algebras over non-archimedean
fields was not studied yet.
\par This article continuous previous investigations of operator
algebras over non-archimedean fields (see \cite{diarlud02} and
references therein), where their spectral theory was described. The
present paper is devoted to investigations of derivations of
operator algebras over infinite spherically complete fields with
non-trivial non-archimedean multiplicative norms having values in
$\Gamma \cup \{ 0 \} $, where $\Gamma $ is a discrete multiplicative
group, particularly over locally compact fields. Theorems about a
differentiation being internal are demonstrated.
\par All results of this paper are obtained for the first time.

\section{Differentiations of operator algebras}
\par {\bf 1. Definitions.} Suppose that $\bf F$ is an infinite field supplied
with a non-archimedean non-trivial multiplicative norm relative to
which it is complete as a uniform space. Let $X$ be a Banach space
over $\bf F$, denote by $L(X)=L(X,X)$ a Banach space of all
continuous $\bf F$-linear operators from $X$ into $X$.
\par An algebra $\Psi $ contained in $L(X)$ over $\bf F$ such that $A^t\in \Psi $ for
each  $A\in \Psi $ will be called an algebra with transposition,
where a mapping $A\mapsto A^t$ on $\Psi $ is called a transposition,
if it is $\bf F$-linear and $(A^t)^t=A$ and $(AB)^t=B^tA^t$ for each
$A, B\in \Psi $.
\par A Banach algebra with transposition is called a $T$-algebra.
\par A subalgebra of $L(X)$ will be called an operator algebra.
\par An operator $A\in \Psi $ is called symmetric if $A^t=A$.
\par An $\bf F$-linear continuous mapping $D: \Psi \to \Psi $ on an algebra $\Psi $ over
$\bf F$ is called a derivation of $\Psi $ if $D(AB)=D(A)B+AD(B)$ for
each $A, B\in \Psi $.
\par For subalgebras $\Phi $ and $\Psi $ of $L(X)$ satisfying the condition
if $A\in \Phi \cup \Psi $ then $A^t\in L(X)$ let $T \{ \Phi ,\Psi \}
$ denote the minimal $T$ subalgebra in $L(X)$ containing $\Phi $ and
$\Psi $.
\par A norm of an operator $A\in L(X)$ is defined as
\par $ \| A \| := \sup_{x\ne 0} \| Ax \| / \| x \| .$
\par A homomorphism $\phi : \Psi \to L(X)$ such that it is $\bf F$-linear,
$\phi (aA+bB)=a \phi (A)+b\phi (B)$, and multiplicative
$\phi (AB)=\phi (A)\phi (B)$ and $\phi (A^t)=[\phi (A)]^t$ for each
$A, B\in \Psi $ and $a, b\in {\bf F}$ is called a representation of
a $T$ algebra on $X$. If additionally $\phi $ is bijective, then
such representation $\phi $ is called faithful.

\par {\bf 2. Remark.} Henceforth, operator $T$-algebras are considered.
Denote by $c_0(\alpha ,{\bf F})$ the Banach space of all mappings
$x: \alpha \to {\bf F}$ satisfying the condition that for each
$\epsilon >0$ the set $ \{ j: ~ j\in \alpha ; ~ |x_j|>\epsilon \} $
is finite, where $c_0(\alpha ,{\bf F})$ is supplied with norm $ \| x
\| := \sup_{j\in \alpha } |x_j|$, $ ~ x_j=x(j)$, $~ \alpha $ is a
set. That is either $\Gamma _{\bf F} := \{ |x|: x\in {\bf
F}\setminus \{ 0 \} \}$ is discrete or $card (\alpha )<\aleph _0$.
\par If $X=c_0(\alpha ,{\bf F})$ and $A\in L(X)$, then a transposed
operator $A^t$ can be defined by the equality $A^t_{j,k}=A_{k,j}$
for each $j, k\in \alpha $, where $Ae_k=\sum_{j\in \alpha }
A_{k,j}e_j$ with $A_{k,j}\in {\bf F}$, $~e_j\in c_0(\alpha ,{\bf
F})$ denotes the basic vector $e_j(k)=\delta _{k,j}$ for each $k\in
\alpha $, $ ~ \delta _{k,j}= 0$ for $j\ne k$, while $\delta
_{j,j}=1$. For $X=c_0(\alpha ,{\bf F})$ this operation $A\mapsto
A^t$ will serve as the transposition if $A$ and $A^t$ are in $L(X)$.
\par If $\bf F$ is a spherically complete non-archimedean field with discrete
multiplicative group $\Gamma _{\bf F}$ or $X$ is finite dimensional
over $\bf F$, then a Banach space $X$ over $\bf F$ is isomorphic
with $c_0(\alpha ,{\bf F})$ for some set $\alpha $ (see Theorems
5.13 and 5.16 \cite{roo}). Henceforward, it is supposed that either
a spherically complete field $\bf F$ is locally compact or it
contains a family $ \{ {\bf G}_{\alpha } : ~ \alpha \in \mu \} $ of
locally compact subfields ${\bf G}_{\alpha }$ such that their union
is dense in $\bf F$, that is $\overline{\bigcup_{\alpha } {\bf
G}_{\alpha }} ={\bf F}$, where $\bar{A}$ denotes the completion of a
subset $A$ relative to the uniformity inherited from $\bf F$.
\par Henceforth, it is supposed that a Banach space $X$ is
isomorphic with $c_0(\alpha ,{\bf F})$. \par Let $C_{\infty
}(\Lambda ,{\bf F})$ denote a Banach algebra of all continuous
functions $f: \Lambda \to {\bf F}$ such that for each $\epsilon >0$
there exists a compact subset $V$ in $\Lambda $ for which $|f(x)|\le
\epsilon $ for every $x\in \Lambda \setminus V$, where $\Lambda $ is
a zero-dimensional locally compact Hausdorff space, while $C(\Lambda
,{\bf F})$ denotes the algebra of all continuous functions $f:
\Lambda \to {\bf F}$. If a Banach algebra $\Psi $ is isomorphic with
$C_{\infty }(\Lambda ,{\bf F})$, then it is called a $C$-algebra.
\par If $\bf F$ is a field with a multiplicative norm and $\Lambda $
is a subset in $\bf F$, a space of all continuous functions $f:
\Lambda \to {\bf F}$ so that for each $\epsilon >0$ a positive
number $0<r<\infty $ exists for which $|f(x)|<\epsilon $ for each
$x\in \Lambda $ with $|x|>r$ is denoted by $C_{\infty }(\Lambda
,{\bf F})$. That is $C_{\infty }(\Lambda ,{\bf F})$ is a space of
continuous functions tending to zero at infinity.
\par Evidently, each $C$-algebra is a $T$-algebra.
\par If a Banach space $X$ is over a spherically complete field $\bf
F$ and $X^*$ is its topological dual Banach space, i.e. of all
continuous $\bf F$-linear functionals $y^*: X\to {\bf F}$, then each
$A\in L(X,Y)$ has an adjoint operator $A^*: Y^*\to X^*$, where $Y$
is a Banach space over $\bf F$, $A^*\in L(Y^*,X^*)$. On the other
hand, $X^*$ is the Banach space over the spherically complete field
$\bf F$ and hence isomorphic with $c_0(\beta ,{\bf F})$ for some set
$\beta $. But each vector $x\in X=c_0(\alpha ,{\bf F})$ gives rise
to a continuous $\bf F$-linear functional $x^*z :=\sum_{j\in \alpha
}x_jz_j$ for each $z\in X$. Therefore, the natural embedding
$X\hookrightarrow X^*$ exists, that is $\alpha \subset \beta $. This
implies, that the operation $L(X)\ni A\mapsto A^*\in L(X^*)$ can be
considered as an extension of $A\mapsto A^t$ from $X$ onto $X^*$ for
each $A\in L(X,X)=L(X)$ (see also Chapter 3 in \cite{roo}).
\par For a Banach space $X$ over a spherically complete field $\bf
F$ each closed linear subspace $Y$ is orthocomplemented in
accordance with Theorems 5.13 and 5.16 \cite{roo}. Therefore, in
such case it is written below for short a projection $\pi _Y: X\to
Y$ instead of an orthoprojection, where $\pi _Y(X)=Y$ (see also
\cite{diarlud02}).

\par {\bf 3. Lemma.} {\it If $\Phi $ is a $C$-algebra over $\bf F$
and $D$ is its differentiation, then $D=0$ on it.}
\par {\bf Proof.} In the space $C_{\infty }(\Lambda ,{\bf F})$
an $\bf F$-linear subspace of simple functions $$f(x)=\sum_{j=1}^n
a_j\chi _{B_j}$$ is dense, where $a_j\in {\bf F}$, $B_j$ is a clopen
subset in $\Lambda $, $ ~ \chi _B$ denotes the characteristic
function of a subset $B$ in $\Lambda $, that is $\chi _B(x)=1$ for
each $x\in B$, while $\chi _B(x)=0$ for any $x\in \Lambda \setminus
B$. Then $D(\chi _B)=D(\chi _B^2) = 2D(\chi _B)\chi _B$,
consequently, $D(\chi _B)(1-2\chi _B)=0$ and hence $D(\chi _B)=0$. A
differentiation $D$ is $\bf F$-linear and continuous, consequently,
$D(f)=0$ for each $f\in C_{\infty }(\Lambda ,{\bf F})$.

\par {\bf 4. Lemma.} {\it If $A\in L(X)$, where $X$ is Banach space
over a locally compact field $\bf F$, and an operator $A$ is such
that $\overline{{\bf F}(A)}$ is a least closed $C$-subalgebra of
$L(X)$ containing $A$, and a spectrum of $\overline{{\bf F}(A)}$ is
contained in a closed ball $B({\bf F},0,\| A \| )$ containing $0$ in
$\bf F$ of radius $ \| A \| $, $~D: \overline{{\bf F}(A)}\to L(X)$
and $\pi _{\overline{{\bf F}(A)}}D: \overline{{\bf F}(A)}\to
\overline{{\bf F}(A)}$ are differentiations, then $\pi
_{\overline{{\bf F}(A)}}DB=0$ for each $B\in \overline{{\bf F}(A)}$,
where $\pi _{\Psi }: L(X)\to \Psi $ denotes an $\bf F$-linear
projection on a closed subalgebra $\Psi $ in $L(X)$.}
\par {\bf Proof.} The field  $\bf F$ is locally compact,
consequently, it is spherically complete. Therefore, the Banach
subspace $\overline{{\bf F}(A)}$ in $L(X)$ is orthocomplemented in
the non-archimedean sense and the continuous $\bf F$-linear
projection $\pi _{\overline{{\bf F}(A)}}$ exists (see Chapter 5
\cite{roo}). Take a closed ball $B({\bf F},0, \| A \| )$ in the
field $\bf F$, where $B(Y,z,r) := \{ x\in Y: ~ \rho (x,z) \le r \} $
denotes a closed ball with center $z$ of radius $0<r$ in a metric
space $Y$ with a metric $\rho $. Since the field $\bf F$ is locally
compact, this ball $B({\bf F},0, \| A \| )$ is compact. So the
$C$-algebra $C(B({\bf F},0, \| A \| ),{\bf F})$ of all continuous
functions $f: B({\bf F},0,\| A \| )\to {\bf F}$ exists. For each
polynomial $P_n(x)$ of degree $n$ on $B({\bf F},0, \| A \| )$ the
corresponding operator $P_n(A)$ is defined, where $A^0=I$ is the
unit operator, $A^nx=A(A^{n-1}x)$ for each $x\in X$. The $\bf
F$-linear space of polynomials is dense in $C(B({\bf F},0, \| A \|
),{\bf F})$ in accordance with Kaplansky's theorem 43.3
\cite{schikb}. By the conditions of this lemma a spectrum of a
$C$-algebra $\overline{{\bf F}(A)}$ is contained in a closed ball
$B({\bf F},0,\| A \| )$. Therefore, $f(A)$ is defined for each
continuous function $f: B({\bf F},0, \| A \| )\to {\bf F}$ and
$\overline{{\bf F}(A)}$ is contained in $C(B({\bf F},0, \| A \|
),{\bf F})$ as the closed subalgebra. Certainly, $\pi
_{\overline{{\bf F}(A)}} DA \in \overline{{\bf F}(A)}$ and $\pi
_{\overline{{\bf F}(A)}} D(AB)= \pi _{\overline{{\bf F}(A)}} (DA)B+
A\pi _{\overline{{\bf F}(A)}} DB$ for each $A, B\in \overline{{\bf
F}(A)}$, hence $\pi _{\overline{{\bf F}(A)}} D$ is the continuous
differentiation on $\overline{{\bf F}(A)}$, since the operators $\pi
_{\overline{{\bf F}(A)}}$ and $D$ are continuous. A closed
subalgebra of a $C$-algebra is a $C$-algebra by Corollary 6.13
\cite{roo}. Therefore, by the preceding lemma the differentiation
$\pi _{\overline{{\bf F}(A)}}$ on $\overline{{\bf F}(A)}$ is
degenerate.

\par {\bf 5. Definition.} Let $\rho : \Psi \to {\bf F}$ be a linear
continuous functional on a $T$-algebra $\Psi $ over $\bf F$. If
$\rho (A^t)=\rho (A)$ for each $A\in \Psi $, then $\rho $ will be
called symmetric. If a symmetric continuous functional $\rho $ is
such that $\rho (I)=1$, then $\rho $ is called a state of $\Psi $. A
state $\rho $ of a $T$-algebra $\Psi $ is definite on a symmetric
operator $A$, when $\rho (A^n)=[\rho (A)]^n$ for every natural
number $n$.
\par A functional $\rho $ is called multiplicative on a $T$-algebra $\Psi $, if $\rho
(AB)=\rho (A)\rho (B)$ for each $A, B\in \Psi $.

\par {\bf 6. Lemma.} {\it Let $A$ be an operator in a Banach
algebra $\Psi $ over a locally compact field $\bf F$ and let a
continuous linear functional $\rho : \overline{{\bf F}(A)} \to \bf
F$ be multiplicative on $A^n$ for each $n\in {\bf N}$. Suppose that
$\rho $ has a $\bf K$-linear extension on $C(B({\bf F},0, \| A \|
),{\bf K})$, where a field ${\bf K}$ is an extension of a field
${\bf F}$. Then $\rho (f(A))=f(\rho (A))$ for each $f\in C(B({\bf
F},0, \| A \| ),{\bf K})$.}
\par {\bf Proof.} From $\rho (A^n)=[\rho (A)]^n$ for every $n\in {\bf N}$ it follows that $\rho $ is
multiplicative on the Banach algebra $\overline{{\bf F}(A)}$
generated by $A$. But $\overline{{\bf F}(A)}$ is the Banach algebra
having an isometric embedding into $C(B({\bf F},0, \| A \| ),{\bf
F})$. For each polynomial $P_m(x)$ on $B({\bf F},0, \| A \| )$ with
values in $\bf K$ due to multiplicativity and $\bf F$-linearity of
$\rho $ one gets $\rho (P_m(A)) = P_m (\rho (A))$. From Kaplansky's
theorem and continuity of $\rho $ it follows that $\rho (f(A)) =
f(\rho (A))$ for each continuous function $f: B({\bf F},0, \| A \|
)\to {\bf K}$.

\par {\bf 7. Lemma.} {\it Let $A$ be a symmetric operator in a
$T$-algebra $\Psi $ over a locally compact field $\bf F$ and let
${\bf K}$ be its extension so that $ ~ \sqrt[m]{x} \in {\bf K}$ for
every $x\in {\bf F}$ and for each natural number $m\ge 2$, then for
a marked natural number $n\ge 2$ there exists $B\in \Psi _{\bf K}$
such that $B^n=A$, where $\Psi _{\bf K}$ is an extension of an
algebra $\Psi $ over ${\bf K}$.}
\par {\bf Proof.} Suppose that a set $\alpha $ is infinite.
A net of projection operators $\pi _{\gamma }$ on finite dimensional
subspaces $c_0(\gamma ,{\bf F})$ in $c_0(\alpha ,{\bf F})$ exists,
where $\gamma $ are finite subsets in $\alpha $ and their family
$\Upsilon $ is partially ordered by inclusion. Then ${\cal B} = \{
\beta : ~ \beta = \alpha \setminus \gamma ; \gamma \in \Upsilon \} $
is a filter base. For each $x\in X=c_0(\alpha ,{\bf F})$ the limit
$\lim_{\cal B} \pi _{\gamma }x=x$ exists. Therefore, $\lim_{\cal B}
\pi _{\gamma } A\pi _{\gamma }x=Ax$ for each $x\in X$. Each operator
$\pi _{\gamma } A \pi _{\gamma }$ is compact from $X$ into $X$. In
view the decomposition theorem of compact operators (see Lemma 1 and
Note 2 of Appendix A in \cite{ludjmsqim05}) it has the decomposition
\par $(1)$ $\pi_{\gamma }A\pi_{\gamma }=C^{-1}_{\gamma }\Lambda _{\gamma
}C_{\gamma }$ over $\bf K$,\\ where $C_{\gamma }$ is an invertible
operator on $c_0(\gamma ,{\bf K})$ and $\Lambda _{\gamma }$ is a
diagonal operator on $c_0(\alpha ,{\bf K})$ relative to its standard
base, moreover, $C_{\gamma }^t=C_{\gamma }^{-1}$ for a symmetric
operator $\pi _{\gamma }A\pi _{\gamma }$. The latter decomposition
automatically encompasses the case of finite $\alpha $ also. Thus
$P_n(\pi_{\gamma }A\pi _{\gamma })$ is correctly defined for each
polynomial $P_n$ on ${\bf F}$ with values in ${\bf K}$ and
$\lim_{\cal B}P_n(\pi _{\gamma }A\pi _{\gamma })x=P_n(A)x$ for every
vector $x\in X$.

\par From the embedding $\overline{{\bf
F}(A)}\hookrightarrow C(B({\bf F},0, \| A \| ),{\bf F})$ (see \S 6)
and the continuity of the function ${\bf F}\ni x\mapsto
\sqrt[n]{x}\in {\bf K}$ it follows that $B=\sqrt[n]{A}\in
\overline{{\bf K}(A)}$. The latter algebra is contained in $\Psi
_{\bf K}$.

\par {\bf 8. Lemma.} {\it Let $D$ be a derivation operator on a
$T$-algebra $\Psi $ over a field $\bf F$, let also $\bf K$ be an
extension of $\bf F$ complete relative to its uniformity. Then $D$
has a continuous $\bf K$-linear extension on $\Psi _{\bf K}$ as a
derivation operator.}
\par {\bf Proof.} A field $\bf K$ contains $\bf F$ as a subfield, a
multiplicative non-archimedean norm on $\bf F$ has a multiplicative
non-archimedean extension on $\bf K$ (see \cite{roo,schikb}). The
completion $\tilde{\bf K}$ of $\bf K$ relative to this norm is a
field complete relative to the norm. A Banach space $X$ over a field
$\bf F$ has an $\bf F$-linear continuous embedding into
$X_{\tilde{\bf K}}$, where $X = c_0(\alpha ,{\bf F})$ and
$X_{\tilde{\bf K}} = c_0(\alpha , \tilde{\bf K})$. \par Then the
closure of the $\tilde{\bf K}$-linear span of $\Psi $ in
$L(X_{\tilde{\bf K}})$ gives $\Psi _{\tilde{\bf K}}$ such that the
embedding $\Psi \hookrightarrow \Psi _{\tilde{\bf K}}$ is continuous
relative to the operator norm. By the condition of this lemma $\bf
K$ is complete relative to its uniformity, hence $\tilde{\bf K}=\bf
K$ and $\Psi _{\bf K}=\Psi _{\tilde{\bf K}}$. \par Put $DbA=bDA$ for
each $b\in \tilde{\bf K}$ and $A\in \Psi $. Therefore,
$D(bAB)=bD(AB)=b(DA)B+bADB$ for any $b\in \tilde{\bf K}$ and $A,
B\in \Psi $. Moreover, $ \| D(bAB) \| \le |b| \max ( \| (DA)B \|, \|
A(DB) \| )$ and $(DbA+tB)\le \max ( |b| \| DA \| , |t| \| DB \| )$
for each $b, t\in \tilde{\bf K}$ and $A, B\in \Psi $, consequently,
$D$ has a continuous $\tilde{\bf K}$-linear extension as a
derivation operator on the $\tilde{\bf K}$-linear span of $\Psi $ in
$L(X_{\tilde{\bf K}})$ and hence on its completion $\Psi
_{\tilde{\bf K}}.$

\par {\bf 9. Lemma.} {\it Suppose that $D$ is a derivation of a $T$-algebra
$\Psi $ over a field $\bf F$ having an extension up to a derivation
on $\Psi _{\bf K}$, where a field ${\bf K}$ is an extension of ${\bf
F}$. Let $\rho $ be a definite state on a symmetric operator $A$ and
$A=B^2$ for some symmetric operator $B\in \overline{{\bf K}(A)}$ and
$\rho $ has an extension as a state on $\Psi _{\bf K}$ and either
$DB\in \overline{{\bf K}(A)}$ or $\rho (B(DB))=\rho (B)\rho (DB)$
and $\rho ((DB)B)=\rho (DB)\rho (B)$. Then $\rho (DA)=0$.}
\par {\bf Proof.} The differentiation operator $D$ is $\bf F$-linear, hence
$DA = D (A - \rho (A)I)$, since $DI=0$. Therefore, without loss of
generality it is sufficient to consider the case $\rho (A)=0$, since
$\rho (A - \rho (A)I)=0$. Consider an operator $B=A^{1/2}$, i.e
$A=B^2$. A state $\rho $ is multiplicative on the $T$-algebra
$\overline{{\bf F}(A)}$ generated by $A$, since $\rho (A^n)=[\rho
(A)]^n$ for each natural number $n$. Thus $\rho (B)=0$ and hence
$\rho (DA)=\rho ((DB)B)+\rho (B(DB))=\rho (B)\rho (DB)+\rho (DB)\rho
(B)=0$, since either $DB\in \overline{{\bf K}(A)}$ or $\rho
(B(DB))=\rho (B)\rho (DB)$ and $\rho ((DB)B)=\rho (DB)\rho (B)$.

\par {\bf 10. Theorem.} {\it Suppose that $D$ is a derivation of a
$T$-algebra $\Psi $ over a locally compact field $\bf F$ such that $
~ \sqrt{x} \in {\bf K}$ for each $x\in {\bf F}$, where a field ${\bf
K}$ is an extension of ${\bf F}$, and that $\cal Z$ is a center of
$\Psi $. Then $D$ annihilates $\cal Z$.}
\par {\bf Proof.} Since $\Psi $ is a Banach algebra, its center $\cal Z$ is
closed in $\Psi $. The field $\bf F$ is spherically complete,
consequently, $\cal Z$ is complemented in $\Psi $. Take $A\in \cal
Z$ and consider $\overline{{\bf F}(A)}$ which is complemented in
$\cal Z$ and hence in $\Psi $. Any element $A\in \Psi $ can be
presented as $A=A_1-A_2$, where $A_1^t=A_1$ is symmetric and
$A_2^t=-A_2$ is antisymmetric, $A_1=\frac{A^t+A}{2}$,
$A_2=\frac{A^t-A}{2}$. If $A\in \cal Z$, then $A_1, A_2\in \cal Z$.
 Then $\sqrt{-1}=i\in {\bf K}$ and
$(iA_2)^t=-A_2$, where $i={{~0~1}\choose{-1~0}}$ over the field $\bf
F$, when $i\notin {\bf F}$. On the other hand, if $i\in {\bf F}$,
one can consider $\Psi \oplus \Psi $ on $X\oplus X$, where
${{~0~~A_2}\choose{-A_2~0}}$ is symmetric, when $A_2$ is
antisymmetric. Therefore, it is sufficient to consider the case of
symmetric $A\in \cal Z$.
\par Choose any multiplicative $\bf F$-linear continuous functional $\rho
$ on $\overline{{\bf F}(A)}$ so that $\rho (I)=1$ and $\rho (A)\ne
0$. Consider a projection $\pi _YDA$ of $DA$ onto a Banach subspace
$Y=\Psi \ominus \overline{{\bf F}(A)}$, i.e. $\Psi =Y\oplus
\overline{{\bf F}(A)}$. Take any continuous extension of $\rho $ so
that $\rho (\pi _YDA)\ne 0$ and such that $\rho $ is multiplicative
on $\overline{{\bf F}(A,\pi _YDA)}$, where $\overline{{\bf
F}(A_1,...,A_n)}$ denotes a minimal closed subalgebra of $\Psi $
containing elements $A_1,...,A_n\in \Psi $. This is possible, since
$\overline{{\bf F}(A,\pi _YDA)}$ is the algebra with two commuting
generators $[A,\pi _YDA]=0$. Moreover, the inclusion $A\in \cal Z$
implies $\overline{{\bf F}(A)} \subset \cal Z$ and $\Psi /{\cal Z}=
(\Psi /\overline{{\bf F}(A)})/({\cal Z}/\overline{{\bf F}(A)})$ and
$(\pi_YDA)+\overline{{\bf F}(A)}=\theta (\pi_YDA) = DA +
\overline{{\bf F}(A)}$, where $\theta : \Psi \to \Psi /
\overline{{\bf F}(A)}$ denotes the quotient mapping. Then we also
get $\theta (C) = C + \overline{{\bf F}(A)}$ and
$D(AB)+\overline{{\bf F}(A)} =(DA)B +A(DB) + \overline{{\bf F}(A)} =
(DA+ \overline{{\bf F}(A)})(B+\overline{{\bf F}(A)})+
(A+\overline{{\bf F}(A)})(DB+\overline{{\bf F}(A)})= \theta (D(AB))=
\theta (DA)\theta (B) + \theta (A)\theta (DB)$, consequently,
$\theta \circ D$ is the differentiation on the quotient algebra
$\Psi /\overline{{\bf F}(A)}$. \par If $V\in \overline{{\bf F}(\pi
_YDA)}\ominus \overline{\bf F}(A)$ is a non zero element and $V^n\in
\overline{\bf F}(A)$ for some natural number $n\ge 2$, then take an
algebraically closed field ${\bf K}$ containing ${\bf F}$ so that
$\sqrt[n]{x}\in {\bf K}$ for each $x\in {\bf F}$. Therefore, $V\in
\overline{\bf K}(A)$ and one can take $\rho (V)=\sqrt[n]{f(\rho
(A))}$ with $Q=V^n=f(A)\in \overline{\bf F}(A)$, where $f$ is a
continuous function from $B({\bf F},0,1)$ into ${\bf F}$ (see Lemmas
6 and 7). Thus it remains to treat the variant when $V^n\notin
\overline{{\bf F}(A)}$ for each natural number $n$.
\par For this it is sufficient to choose a multiplicative extension
of $\rho $ on $\overline{{\bf F}(\pi _YDA)}\ominus \overline{\bf
F}(A)$ putting $\rho (V^n)=[\rho (V)]^n\ne 0$ for each natural
number $n$ and for some non zero element $V\in \overline{{\bf F}(\pi
_YDA)}\ominus \overline{\bf F}(A)$. Indeed, without loss of
generality using multiplication on scalars $A\mapsto bA$ for $b\in
{\bf F}\setminus \{ 0 \} $ it is possible to restrict on the case
$\max ( \| A \| , \| V \| )<1$ and choose $0<|\rho (A) | \le \| A \|
$ and $0< |\rho (V)| \le \| V \| $. The Banach subspace
$\overline{{\bf F}(A,\pi _YDA)}$ is closed in $\Psi $, consequently,
by the non-archimedean Hahn-Banach theorem over the spherically
complete field $\bf F$ a functional $\rho $ has a continuous
extension on $\Psi $ (see \cite{roo} and \S 8.203 in \cite{nari}).
\par The family of such functionals $\rho $ separates different elements
of $\Psi $ and $Q\in \Psi \ominus \overline{{\bf F}(A,\pi _YDA)}$,
hence $\rho (DA)=0$ for each such $\rho $ if and only if $DA=0$.
\par Applying Lemmas 6-9 we get the statement of this theorem.

\par {\bf 11. Definition.} If $\Psi $ is a $T$-algebra on a Banach
space $X$ over a field $\bf F$, its strong topology is characterized
by a base of neighborhoods $V_{x_1,...,x_n;\epsilon } := \{ A\in
\Psi : ~ \| Ax_j \| <\epsilon ~ \forall j=1,...,n \} $ of zero,
where $x_1,...x_n\in X$, $~\epsilon >0$, $~ n\in {\bf N}$. If a
field $\bf F$ is spherically complete and $X^*$ is a topological
dual space of $X$, a weak topology on $\Psi $ is given by a base of
neighborhoods $W_{x_1,...,x_n;y_1,...,y_n;\epsilon } := \{ A\in \Psi
: ~ |y_j Ax_j | <\epsilon ~ \forall j=1,...,n \} $ of zero, where
$x_1,...x_n\in X$, $~y_1,...,y_n\in X^*$, $~\epsilon >0$, $~ n\in
{\bf N}$. Denote by $\bar \Psi $ the completion of $\Psi $ relative
to the weak topology.

\par {\bf 12. Lemma.} {\it Suppose that $D$ is a derivation of
a $T$-algebra $\Psi $ on a Banach space $X$ over a spherically
complete field $\bf F$. Then a unique weakly continuous extension
$\bar{D}$ of $D$ on $\bar \Psi $ exists.}
\par {\bf Proof.} The mappings $A\mapsto \frac{A^t+A}{2}=:A_1$ and
$A\mapsto \frac{A^t-A}{2}=:A_2$ are continuous on a $T$-algebra
$\Psi $. An extension $\bf K$ from Lemma 7 of a spherically complete
field $\bf F$ can be considered as an $\bf F$-linear space. By Lemma
8 $D$ has a continuous extension on $\Psi _{\bf K}$ as a derivation
operator. As in Lemma 10 it is sufficient to consider a symmetric
operator $A$. \par Put ${\bf S} := \{ A\in \Psi : ~ \| A \| \le 1, ~
A^t=A \} $ to be the unit ball of symmetric operators. Then the
mapping ${\bf S}\ni A\mapsto y(D(A^2)x)=y(ADAx+(DA)Ax)$ is strongly
continuous at zero, since $|y(ADAx+(DA)Ax)|\le \| D \| \max ( \| Ax
\| \| y \| ; \| x \| \| Ay \| ),$ where $x, y \in X$ and $X$ is
embedded into $X^*$. On the other side, the mapping ${\bf S}\ni
A\mapsto A^{1/2}$ is strongly continuous at zero, since $ \|
A^{1/2}x \| \le \| Ax \| \| x \| $ due to Formula 7$(1)$, where
$x\in X$. Thus $A\mapsto y(DA_1x)-y(DA_2x)=y(DAx)$ is strongly
continuous at zero on $\bf S$. This implies that $H := {\bf S}\cap
q^{-1}(B({\bf F},0,r))$ is strongly closed in ${\bf S}$, where
$q(A):=y(DAx)$ for some marked vectors $x, y\in X$, $~ 0<r<\infty $.
\par Recall that a subset $U$ of a topological $\bf F$-linear space
$Q$ is called absolutely $\bf F$-convex if $B({\bf F},0,1)U+B({\bf
F},0,1)U\subset U$. \par The norm on $\bf F$ is non-archimedean,
i.e. $|a+b|\le \max (|a|,|b|)$ for each $a, b\in {\bf F}$. It can be
lightly seen, that the set $H$ is absolutely $\bf F$-convex and
strongly closed, consequently, $H$ is weakly closed in $\bf S$.
Indeed, if a net $T_n\in H$ strongly converges to $T\in H$, then
$T_n-T\in H$ for each $n$ and hence the net $(T_n-T)$ strongly
converges to zero. Therefore, $y((T_n-T)x)$ converges to zero for
each $x\in X$ and $y\in X^*$.
\par By the non-archimedean Hahn-Banach theorem 8.203 \cite{nari} the set $H$ is closed relative
to a weak topology with functionals from $X$, since the set of
continuous $\bf F$-linear functionals $y\in X$ separates points in
$X$, i.e. from $\lim_n y((T_n-T)x)=0$ for each $y\in X$ it follows
$\lim_n y((T_n-T)x)=0$ for every $y\in X^*$, where $T, T_n\in H$.
\par Therefore, the derivation $D$ is weakly continuous on $B(\Phi
,0,1)$, since the mapping $B(\Phi ,0,1)\ni A\mapsto y(DAx)$ is
continuous for each marked $x, y \in X$ and the derivation operator
$D$ is $\bf F$-linear. This means that $D$ is uniformly continuous
relative to the weak uniformity on $B(\Phi ,0,1)$ and implies that
$D$ has a continuous extension on $\overline{B(\Phi ,0,1)}$ and
hence on $\bar{\Phi }$ with range in $\bar \Phi $ by Theorem 8.3.10
\cite{eng}, since $\overline{B(\Phi ,0,1)}$ is the closed absorbing
set in $\bar{\Phi }$.
\par This extension is $\bf F$-linear as well, since
\par $lim_n D(bT_n+H_n)=b\lim_n DT_n+\lim_n DH_n$ for each
$b\in \bf F$ and $T_n, H_n\in \Phi $ with $\lim_n T_n=T\in \bar{\Phi
}$ and $\lim_nH_n=H\in \bar{\Phi }$. Moreover, $y(D(T_nH_n)x)=
y((DT_n)H_nx)+y(T_n(DH_n)x)$ for each $x\in X$ and $y\in X^*$,
consequently,
\par $ \lim_n \lim_k y((DT_n)H_kx)+y(T_n(DH_k)x)- y(D(T_nH_k)x) =
y((DT)Hx)+ y(T(DH)x)-y(D(TH)x)=0 $, since $D$ is the bounded
operator on $\Phi $ and weakly continuous on $\bar{\Phi }$, hence
$D$ is the derivation on $\bar{\Phi }$ as well.
\par {\bf 13. Definitions.} A $T$-algebra of bounded operators on a
Banach space $X$ over a spherically complete field $\bf F$ closed
relative to the weak operator topology and containing the unit
operator will be called a $W^t$-algebra. For an operator $A\in L(X)$
and a $W^t$-algebra $\Psi $ let $\overline{co}_{\Psi }(A)$ denote
the closure relative to the weak operator topology of finite
combinations $b_1B_1+...+b_nB_N$ of operators $B_j=V_jAV_j^t$, where
$V_j$ is an isometry operator on $X$ for every $j$, i.e. $\| V_jx \|
= \| x \| $ for each $x\in X$, $~b_1,...,b_n\in B({\bf F},0,1)$.
\par If $\Upsilon $ is a family of operators in $L(X)$, then
$\Upsilon ' := \{ C: ~ C\in L(X); ~ [C,T]=0 ~ \forall T\in \Upsilon
 \} $ denotes the commutant of $\Upsilon $, where $[C,T]=CT-TC$
 is the commutator of two operators.
\par A center $Z(\Psi )$ of an algebra $\Psi $ is a set of all its elements
commuting with each element in $\Psi $. An element $A\in Z(\Psi )$
in the center is called central.

\par {\bf 14. Lemma.} {\it Let $A$ be a linear continuous operator $A: X\to X$
on a Banach space over a spherically complete field $\bf F$ and let
$\bf G$ be a locally compact field contained in $\bf F$. Suppose
that $f\in C_{\infty }({\bf F},{\bf F})$ is a continuous function
tending to zero at infinity the restriction of which $f|_{\bf G}$
belongs to $C_{\infty }({\bf G},{\bf G})$. Then a linear continuous
bounded operator $f(A)\in L(X)$ exists.}
\par {\bf Proof.} The field ${\bf Q}_p$ of $p$-adic numbers is locally compact.
Let ${\bf G}$ be a locally compact field so that ${\bf Q}_p\subset
{\bf G}\subset {\bf F}$, i.e. either a locally compact field ${\bf
G}$ containing ${\bf Q}_p$ or the $p$-adic field itself. Then each
$\bf F$-linear operator is also ${\bf G}$ linear. \par Let $X_{\bf
G}$ denote the Banach space over ${\bf G}$ obtained from the Banach
space $X$ over ${\bf F}$ considering ${\bf F}$ as the Banach space
over ${\bf G}$, i.e. by the restriction of the field of scalars.
Take $P$ a projection $P\in {\sf P}_{\bf G}$ on a finite-dimensional
over ${\bf G}$ subspace in $X_{\bf G}$ with ${\sf P}_{\bf G}$
denoting the family of all projections having finite dimensional
ranges partially ordered by inclusion of their ranges in $X_{\bf
G}$. Then each operator $PAP$ can be reduced to the diagonal form
\par $(1)$ $PAP= CTE$ \\ over ${\bf G}$ by a lower and upper
triangular operators $C$ and $E$ respectively invertible on $PX$
with diagonal operator $T$ such that $(C-I)$ and $(E-I)$ are
nilpotent operators on $PX_{\bf G}$ (see Lemma 1 of Appendix A in
\cite{ludjmsqim05}).
\par In accordance with E. Zermelo's theorem on each set $\Lambda $ a relation
exists, which well orders $\Lambda $ (see \cite{eng}). Suppose that
$P_{\beta }$ is a family of projections on a Banach space over a
spherically complete field $\bf F$, where $\beta \in \Lambda $ and a
set $\Lambda $ is well ordered and $P_{\alpha }\le P_{\beta }$ for
each $\alpha \le \beta $. Denote by $\wedge_{\alpha \in \Lambda }
P_{\alpha }$ an projection from $X$ onto the subspace
$\bigcap_{\alpha \in \Lambda } P_{\alpha }X$, while defining
$\vee_{\alpha \in \Lambda } P_{\alpha }:= I- \wedge_{\alpha \in
\Lambda } (I-P_{\alpha })$, where $I$ is the unit operator on $X$,
$~Ix=x$ for each $x\in X$. Then the family $Q_{\alpha } := P_{\alpha
} - \vee_{\beta <\alpha } P_{\beta }$ consists of mutually
orthogonal projections on $X$ such that its sum is $\vee_{\beta \in
\Lambda } Q_{\beta } = \vee_{\beta \in \Lambda } P_{\beta }=:P$.
\par Indeed, $Q_{\beta }\perp P_{\alpha }$ are orthogonal for each $\alpha <\beta $
and $Q_{\beta }\perp Q_{\alpha }$, since $Q_{\alpha }\subseteq
P_{\alpha }$, i.e. $Q_{\alpha }X\subseteq P_{\alpha }X$. Therefore,
$\vee_{\alpha \in \Lambda } Q_{\alpha } \subseteq P$. If $\alpha _1$
is the least element of $\Lambda $, then $P_{\alpha _1}=Q_{\alpha
_1}$. Suppose that \par $P_{\beta } = \vee_{\alpha \le \beta , ~
\alpha \in \Lambda } Q_{\alpha }$ \\ for each $\beta < \gamma \in
\Lambda $. From the definition of $Q_{\gamma }$ it follows, that
\par $Q_{\gamma } = P_{\gamma } - \vee_{\alpha <\gamma , ~
\alpha \in \Lambda } Q_{\alpha }$, consequently, \par $P_{\gamma } =
I-\wedge_{\alpha \le \gamma; \alpha \in \Lambda } (I-Q_{\alpha }) =
\vee_{\alpha \le \gamma , ~ \alpha \in \Lambda } Q_{\alpha }$. \\
Thus by transfinite induction the latter equality is fulfilled for
each $\gamma \in \Lambda $, hence $P\subseteq \vee_{\alpha \in
\Lambda }Q_{\alpha }$, together with the opposite inclusion this
implies $P=\vee_{\alpha \in \Lambda }Q_{\alpha }$.
\par The field ${\bf F}$ is spherically complete and considered as a
Banach space over ${\bf G}$ is isomorphic with $c_0(\beta ,{\bf G})$
for some set $\beta $ by Theorems 5.13 and 5.16 \cite{roo}.
Therefore, $\lim_{{\cal P}_{\bf G}} PAPx = Ax$ for each $x\in X$. To
an operator $Y\in L(X)$ an operator $Y_{\bf G}\in L(X_{\bf G})$
corresponds such that to each matrix element $e_j^*Ye_k$ over $\bf
F$ an operator block on $c_0(\beta ,{\bf G})$ is posed.
\par Then $C-I$ and $E-I$ are nilpotent operators such that
$(C-I)^l=0$ and $(E-I)^l=0$ for each $l\ge m$, where $m$ is an order
of a square $m\times m$ matrix with entries in $\bf G$, i.e.
$m=dim_{\bf G}PX_{\bf G}$ is a dimension of $PX_{\bf G}$ over the
field $\bf G$, operators $C$ and $E$ are as in Formula 14$(1)$.
Therefore,
$$(2)\quad C^k=\sum_{0\le h\le \min (m,k)} {k\choose{h}}(C-I)^h,$$
where $(C-I)^0=I$ is the unit operator, as usually ${k\choose{h}} =
k!/(h!(k-h)!)$ denotes the binomial coefficient. Since
${k\choose{h}}$ are integers, it follows that $|{k\choose{h}}|_{\bf
G} \le 1$ and hence $ \| S(C) \| \le \sup_{0\le h\le \min(m,n)}
|s_h| \| C-I \|^h<\infty $ for each polynomial $$(3)\quad S(x) =
\sum_{k=0}^n s_k x^k$$ on $\bf G$ with coefficients $s_k\in \bf G$,
$~ s_n\ne 0$, of degree $n=deg ~ S$. Moreover, $S(T)= diag
(S(t_1),...,S(t_m))$ for the diagonal operator $T=diag
(t_1,...,t_n)$ in the corresponding non-archimedean orthonormal
basis in the subspace $PX_{\bf G}$ over the field $\bf G$, where
$t_1,...,t_m\in \bf G$. On the other hand, applying Theorems 5.4,
5.11 and 5.16 \cite{diarlud02} we get:
$$(4)\quad \| S (PAP) \| \le \sup_{t\in {\bf G}, |t|\le \|
PAP \| } | S(t) | ,$$ since $ \| PAP \| = \sup_{1\le v, l \le m}
|q_v^* PAP q_l|$, $ ~ \| P \| =1$ for each non-degenerate projection
operator, where $q_j$ is a non-archimedean orthonormal basis in
$PX_{\bf G}$, $~q_j^*\in P{X_{\bf G}}'$ denotes a $\bf G$ linear
functional corresponding to $q_j$.
\par In view of Kaplansky's theorem a family of polynomials
is dense in $C(B({\bf G},0,r),{\bf G})$ for each $0<r<\infty $ for
the locally compact field, since the ball $B({\bf G},0,r)$ is
compact. For every $f\in C_{\infty }({\bf G},{\bf G})$ and each
$r=p^j\in \Gamma _{\bf G} := \{ |x|: ~x\in {\bf G}\setminus \{ 0 \}
\} $ a sequence $\{ S_{n_j(k)}: ~ k \} $ of polynomials exists
uniformly converging to $f$ on $B({\bf G},0,p^j)$, where
$n_j(k)<n_j(k+1)$ for each $k\in {\bf N}$, $~n=deg ~ S_n$. By
induction construct them such that $\{ n_{j+1}(k): ~ k\in {\bf N} \}
\subset \{ n_j(k): ~ k\in {\bf N} \} $ for each natural number $j\in
{\bf N}$. Choosing the diagonal subsequence $\{ n_j(j): ~ j \in {\bf
N} \} $ one gets a sequence of polynomials $S_{n_j(j)}$ point wise
converging to $f$ on $\bf G$ and uniformly converging to $f$ on each
bounded ball $B({\bf G},0,r)$, since $\lim_{|t|\to \infty } f(t)=0$.
Since $ \| A \| <\infty $, the function $$(5)\quad f(A)x =
\lim_{j\to \infty } \lim_{P\in {\sf P}_{\bf G}} S_{n_j(j)} (C)
S_{n_j(j)} (T) S_{n_j(j)}(E)x$$ exists for each $x\in X$, where $C$,
$T$ and $E$ correspond to $PAP$, $~P\in {\sf P}_{\bf G}$. Evidently
it is linear by $x\in X$, since $\lim_{P\in {\sf P}_{\bf G}}
S_{n_j(j)} (C)S_{n_j(j)} (T)S_{n_j(j)}(E)$ is a linear operator on
$X$ over $\bf F$ for each $j$. Since ${\bf G}\subset {\bf F}$,
Formulas $(1-5)$ imply that
$$ (6)\quad \| f(A) \| \le \sup_{t\in {\bf F}, ~ |t|\le \| A \| } |f(t)|< \infty .$$

\par {\bf 15. Theorem.} {\it Suppose that $\Phi $ is an algebra with transposition of
bounded linear operators on a Banach space $X$ over a spherically
complete field $\bf F$, then each $A\in B(\bar{\Psi },0,1)$ belongs
to the strong operator closure $\overline{B(\Psi ,0,1)}$ of the unit
ball $B(\Psi ,0,1)$ of $\Psi $. If $Q$ is a symmetric operator in
$B(\bar{\Psi },0,1)$, then $Q$ is in the strong-operator closure of
the set of symmetric operators in $B(\Psi ,0,1)$.}
\par {\bf Proof.} As in section 12 for an absolutely convex subset
$E$ of $L(X)$ the weak- and strong-operator closures coincide, since
$X$ is a Banach space over a spherically complete field $\bf F$.
Indeed, for each proper norm closed linear subspace $Y$ of $X$ and a
point $x\in X\setminus Y$ a continuous linear functional $f: X\to
\bf F$ exists such that $f(x)=1$ and $f(Y)=0$ due to the Hahn-Banach
theorem over $\bf F$ (see \S 8.203(f) \cite{nari}). For each point
$x$ outside the norm closure $cl_n U$ of a subset $U$ in $X$ there
exists a closed ball $B(X,x,r) := \{ z\in X: ~ \| z-x \| \le r \} $
containing $x$ of radius $0<r<\infty $ such that the intersection
$(cl_n U)\cap B(X,x,r) = \emptyset $ is void with $r\in \Gamma _{\bf
F}$, where $\Gamma _{\bf F} := \{ |b|: ~b \in {\bf F}\setminus \{ 0
\} \} $ is a multiplicative group contained in $\bf R$. The
multiplicative norm on $\bf F$ is non-trivial, consequently, zero is
the limit point of $\Gamma _{\bf F}$ in $\bf R$. Therefore, a radius
$r>0$ can be chosen so that $\inf_{y\in cl_n U} \| x-y \|
>r$. \par If $V$ is an absolutely convex norm closed subset of $X$ and $x\in X\setminus V$,
there exists a hyperplane $y+Y$ in $X$ which does not contain $x$
and does not intersect $V$, where $~ y = \lambda x$ for some
$\lambda \in \bf F$, $0< |\lambda | \le 1$, $ ~ X=Y\oplus {\bf F}$.
The topological dual space $X'$ of all continuous linear functionals
$f: X\to \bf F$ separates points in $X$, consequently, there exists
a family $ \{ f_{\beta } \} \subset X'$ of continuous linear
functionals and closed subsets $K_{\beta }$ in the field $\bf F$
such that $V = \bigcap_{\beta } f^{-1}_{\beta } (K_{\beta })$.
\par Evidently, if $A$ is in the strong operator closure of $E$,
then it is in the weak operator closure of $E$. Let now $A$ be in
the weak-operator closure of $E$. Consider vectors $x_1,...,x_n\in
X$ and the $n$-fold direct sum $X^{\oplus n} = X\oplus ... \oplus
X$. An operator $G$ on $X$ induces and operator ${\tilde G}=G\oplus
... \oplus G$ on $X^{\oplus n}$. Therefore, $\{ \tilde{G}: G\in E \}
=: \tilde{E}$ is an absolutely convex subset of $X^{\oplus n}$,
hence $\tilde{E}\tilde{x}$ is an absolutely convex subset of
$X^{\oplus n}$, where $\tilde{x}=(x_1,...,x_n)$. If $\tilde{A}$ is
in the weak-operator closure of $\tilde E$, $\tilde{A}\tilde x$ is
in the weak closure of $\tilde{E}\tilde x$, hence in the norm
closure of $\tilde{E}\tilde x$ in $X^{\oplus n}$ due to the fact
demonstrated above. \par That is for each $\epsilon >0$ there exist
$T\in E$ such that $ \| Tx_j-Ax_j \| <\epsilon $ for each
$j=1,...,n$. Thus the weak-operator closure and the strong-operator
closure of $E$ coincide.
\par In view of Lemma 14, an operator $f(A)$ is defined for each symmetric
bounded operator $A\in L(X)$ and hence $f(A)$ for each $A$ in
$\bar{\Psi }$, since $\lim_{|t|\to \infty } f(t)=0$. Moreover, for
each bounded symmetric operator $A$ a symmetric operator $A_{\bf G}$
on $X_{\bf G}$ corresponds, since $x^t=x$ for each $x\in {\bf F}$.
\par Let $Q$ be a symmetric operator in $\bar{\Psi }$, let also
$K_b$ be a net of operators in $\Psi $ weak-operator converging to
$Q$. Then $(K_b+K_b^t)/2$ is a net of symmetric operators in $\Psi $
converging to $Q$ relative to the weak-operator topology. But the
set of symmetric operators in $\Psi $ is absolutely convex and from
the fact demonstrated above $Q$ is in its strong-operator closure.
\par Consider a symmetric operator $Q\in B(\bar{\Psi },0,1)$ and a net of
symmetric operators $M_b\in \Psi $ strong-operator converging to
$Q$. Let $p$ be a prime number so that ${\bf F}$ is an extension of
the $p$-adic field ${\bf Q}_p$, hence up to an equivalence of
multiplicative norms on $\bf F$ we have $|p| := |p|_{\bf F} = 1/p$
(see \cite{schikb,wei}). Take a continuous function $f: {\bf F}\to
{\bf F}$ so that $f(t)=t$ on $B({\bf F},0,1)$, while
$f(t)=p^{2k-1}t$ on $B({\bf F},0,p^k)\setminus B({\bf F},0,p^{k-1})$
for each natural number $k\in {\bf N} := \{ 1, 2, 3,... \} $, since
the ball $B({\bf F},0,r)$ is clopen (simultaneously closed and open)
in $\bf F$, where $r>0$. The function $f$ has the natural extension
on the field $\bf K$ containing $\bf F$ so that $\sqrt[n]{x}\in \bf
K$ for each $x\in \bf F$, putting $f(t)=t$ on $B({\bf K},0,1)$,
while $f(t)=p^{2k-1}t$ on $B({\bf K},0,p^k)\setminus B({\bf
K},0,p^{k-1})$ for every $k\in {\bf N}$. Since $sp (Q)\subset B({\bf
K},0,1)$ (see \cite{roo,diarlud02}), it follows that $f(Q)=Q$.
Moreover, the function $f$ is strong-operator continuous on the set
of symmetric operators in $\bar{\Psi }$. The inequality $ |f(t)|\le
1$ for each $t$ implies that $ \| f(M_b) \| \le 1$ for each $b$. If
$x, y\in {\bf K}$ and $|x-y|< |x|$, then $|y| = |x|$ due to the
non-archimedean inequality $|x+y|\le \max (|x|, |y|)$ for each $x, y
\in {\bf K}$. Therefore, $f(M_b)$ is strong-operator converging to
$f(Q)$, since $\lim_b S_{n_j(j)}(M_b)= S_{n_j(j)}(Q)$ for each $j$
and $P\in {\sf P}_{\bf G}$. Thus $Q$ is in the strong-operator
closure of the set of symmetric operators from $B(cl_n \Psi ,0,1)$
and hence the strong operator limit of symmetric elements in $B(\Psi
,0,1)$.
\par Generally if $A\in B(\bar{\Psi },0,1)$, then form an operator $A' :=
{{0 ~~ A}\choose{A^t~0}}$ on $X\oplus X$ which is symmetric. Then
$A'\in B(\bar{\Psi _2},0,1)$, where $\Psi _2$ denotes the family of
all operators on $X\oplus X$ presented as $2\times 2$ matrices with
entries in $\Psi $. From the proof above it follows that $A'$ is in
the strong-operator closure of $\Psi _2$. Particularly each entry of
$A'$ is in the strong-operator closure of $B(\Psi ,0,1)$, since each
entry in $B(\Psi _2,0,1)$ is in $B(\Psi ,0,1)$.

\par {\bf 16. Definition.} A derivation $D$ of a subalgebra $\Upsilon $ in
$L(X)$ is called spatial, if an operator $B\in L(X)$ exists such
that $D=ad ~ B|_{\Upsilon }$.

\par {\bf 17. Theorem.} {\it Let $\Psi $ be a $T$-algebra on a
Banach space over a spherically complete field $\bf F$, let also $D$
be a derivation of $\Psi $. Then for each commutative
$W^t$-subalgebra $\Phi $ in a commutant $\Psi '$ a bounded $\bf
F$-linear operator $B=B_{\Phi }\in L(X)$ exists such that $B$
commutes with $\Phi $ and $D=ad ~ B|_{\Psi }$.}
\par {\bf Proof.} Evidently $\Upsilon '$ from Definition 13 is weakly
closed in $L(X)$, particularly, $\Psi '$ is weakly closed. Let $\Xi
$ be a maximal commutative subalgebra of $\Psi '$, hence it is
weakly closed in the commutant $\Psi '$. Consider a lattice $\cal P$
of projection operators in $\Xi $ which corresponds to $\Psi $ (see
Theorems 5.4, 5.11 and 5.16 in \cite{diarlud02}). \par The central
carrier of an operator $A\in \Psi $ is defined to be $(I-P)$, where
$P=\bigcup_{\beta } P_{\beta }$ and $P_{\beta }$ is from the set of
all central projections in $\Psi $ such that $P_{\beta }A=0$, i.e.
every $P_{\beta }$ is in the center $Z(\Psi )$ of $\Psi $. Denote by
$C_A$ the central carrier of $A$, then $C_AA=A$, since $A$ is
continuous and $Ax$ is orthogonal to the range of $P_{\beta }$ for
each $\beta $, but $Range (P_{\beta })\subset Range (P)$.
\par Suppose that $B_{j,k}\in \Psi $ and $Q_{j,k}\in \Psi '$ are
operators, then \par $(i)$ $\sum_k B_{j,k}Q_{k,l}=0$ if and only if
central operators $A_{j,k}\in \Psi $ exist satisfying the
properties: \par $(ii)$ $\sum_k B_{j,k}A_{k,l}=0$ and $\sum_k
A_{j,k}Q_{k,l}=Q_{j,l}$ for each $j, l= 1,...,n$. Particularly,
$BQ=0$ for $B\in \Psi $ and $Q\in \Psi '$ if and only if $C_BC_Q=0$.
\par Indeed, from the properties $$\sum_k B_{j,k}A_{k,l}=0\mbox{ and }\sum_k
A_{j,k}Q_{k,l}=Q_{j,l}$$ of central operators $A_{j,k}\in \Psi $ it
follows that $$\sum_kB_{j,k}Q_{k,j}=\sum_kB_{j,k} \sum_t A_{k,t}
B_{t,j} = \sum_t \sum_k B_{j,k} A_{k,t} B_{t,j} =0.$$
\par On the other side, if $\sum_k B_{j,k}Q_{k,l}=0$, then
one can consider the ring $Mat_n(\Psi ')$ of all $n\times n$
matrices with entries in $\Psi '$ and the union of all projections
$T_n =(A_{j,k})$ in $Mat_n(\Psi ')$ which are annihilated under the
left multiplication $BT_n=0$ by $B$, where $A_{j,k}\in \Psi '$ for
each $j, k$. Consider a diagonal matrix $E_n$ with entries being
projections in $\Psi '$. Then $BE_nT_n=0$, consequently,
$T_nE_nT_n=E_nT_n$ and hence $T_nE_n=(T_nE_nT_n)^t=E_nT_n$. Thus
$A_{j,k}\in Z(\Psi ')$. Then the equality $BQ=0$ implies $T_nQ=Q$,
that is, $\sum_k A_{j,k} Q_{k,l} =Q_{j,l}$ for each $j, l=1,...,n$.
Particularly, if $C_BC_Q=0$, then $BQ=BC_BC_QQ=0$. When $BQ=0$, a
central projection $P$ in $\Psi $ exists such that $PB=0$ and
$PQ=Q$, consequently, $PC_B=0$ and $Range (C_B)\subset Range (P)$,
hence $C_BC_Q=0$.
\par Recall that vectors $y_1,...,y_n,...$ are called mutually orthogonal in the non-archimedean
sense, if $\| t_1y_1+...+t_ky_k \| = \max_{j=1}^k \| t_j y_j \| $
for each $t_1,...,t_k\in {\bf F}$ and $k\in {\bf N}$. Two subspaces
$U$ and $W$ in a normed space $Y$ are called orthogonal and denoted
$U\perp W$ if each vector $x\in U$ is orthogonal to every vector
$y\in W$, $ ~ x\perp y$.
\par A closed $\bf F$-linear subspace $U$ in a normed space $Y$ is
complemented, if a closed $\bf F$-linear subspace $V$ in $Y$ exists
so that $U\cap V = \{ 0  \} $ and $U+V=Y$. It is orthocomplemented
if it is complemented and in addition orthogonal $U\perp V$ to its
complement $V$.
\par We say, that $E_1,...,E_j$ are (mutually) complemented, if $E_lE_k=0$ for each $1\le l\ne k\le j$.
\par A projection operator $E: Y\to Y$ is called an orthoprojection if
$E(Y)\perp E^{-1}(0)$.
\par By Theorem 3.9 \cite{roo} a closed linear subspace $U$ of a Banach space $Y$
is complemented if and only if a projection $P: Y\to U$ exists.
Theorem 3.10 \cite{roo} asserts, that a closed linear subspace $U$
of a Banach space $Y$ over a non-archimedean field is
orthocomplemented if and only if an orthoprojection $E$ of $Y$ on
$U$ exists. In view of Theorems 5.13 and 5.16 \cite{roo} each closed
linear subspace of a Banach space over a spherically complete field
is orthocomplemented. On the other hand, each closed linear subspace
of a Banach space over a spherically complete field has an
orthogonal basis which can be extended to an orthogonal basis of the
entire Banach space. Therefore, without loss of generality we
consider the family $\cal P$ of all orthoprojections $E: X\to X$ (
for short of projections).
\par Then we define a new operator $D_1$ by the formula:
\par $(iii)$ $D_1(A_1E_1+...+A_nE_n) = \bar{D}(A_1)E_1+...+\bar{D}(A_n)E_n$,\\
where $E_1,...,E_n\in \cal P$, $ ~ A_1,...,A_n\in \bar{\Psi }$, $~ n
\in {\bf N}$, $~\bar{D}$ is an extension of $D$ from $\Psi $ onto
$\bar{\Psi }$ in accordance with Lemma 12. If $A_1E_1+..+A_nE_n=0$,
then from the proof above it follows that central operators
$C_{j,k}\in Z(\bar{\Psi })$ exist so that $\sum_{k=1}^n
C_{j,k}E_k=E_j$ and $\sum_{j=1}^n A_jC_{j,k}=0$. In view of Theorem
10 $\sum_{j=1}^n \bar{D}(A_j) C_{j,k} =0$, consequently, $\sum_j
\bar{D}(A_j)E_j=0$ by $(i,ii)$. This means that $D_1$ is
single-valued. Denote by $\Phi $ an algebra over $\bf F$ of all
elements of the form $A_1E_1+...+A_nE_n$ with $A_j$ and $E_j$ as
above. It is indeed an algebra, since $A_jE_jA_kE_k=A_jA_kE_jE_k$
for each $j, k$.
\par The definition of $D_1$ implies that this operator is $\bf F$-linear
and bounded on $\Phi $ due to Formula $(iii)$.  Next we verify, that
$D_1$ is a derivation of $\Phi $.
\par If projections $E_1,...,E_j$ are complemented, take
$F_{j+1} = E_{j+1} -E_{j+1}(E_1+...+E_j)$ and so on by induction.
From $E_{l}(X)\perp E_{l}^{-1}(0)$ for each $l=1,...,j+1$ and
$F_{j+1}= (I-E_1-...-E_j) E_{j+1}$  it follows, that
$(I-E_1-...-E_j)(X)\perp (I-E_1-...-E_j)^{-1}(0)$ and
$(I-E_1-...-E_j)(E_{j+1}X)\perp
E_{j+1}^{-1}(I-E_1-...-E_j)^{-1}(0)$, consequently, $F_{j+1}$ is
also the projection. Then $A_lE_l+A_{j+1}E_{j+1} =
A_l(E_l-E_{j+1}E_l)+(A_l+A_{j+1})E_{j+1}E_l+
A_{j+1}(E_{j+1}-E_{j+1}E_j)$ for each $l\le j$ by induction,
consequently, this induces the decomposition
$A_1E_1+...+A_nE_n=B_1F_1+...+B_nF_n$ with complemented projections
$F_1,...,F_n\in \cal P$ and $B_1,...,B_n\in \bar{\Psi }$.

\par When $E_1,...,E_n$ are complemented projections and $x=\sum_{j=1}^nE_jx$
is a vector in $X$ of unit norm $ \| x \| =1$, then $$ \|
(A_1E_1+...+A_nE_n)x \| = \max_{j=1}^n \| A_jE_j x \| ,$$ since
$A_jE_jx=E_jA_jx$ are mutually orthogonal in the non-archimedean
sense vectors. Moreover, \par $\| A_jE_j x \| \le \| A_jE_j  \| \|
E_j x \| \le \max_{l=1}^n \| A_lE_l  \| $, \\ since $ \max_j \| E_j
x_j \| = \| x \| =1$, hence \par $ \| A_1E_1+...+A_nE_n \| \le
\max_{j=1}^n \| A_jE_j \| $. At the same time
\par $\max_j \| A_j E_j \| \le \| A_1E_1+...+A_nE_n \| $ \\
due to the non-archimedean orthogonality of $E_j$. This implies
\par $ \| D_1 (A_1E_1+...+A_nE_n) \| = \max_{j=1}^n
\| (\bar{D}A_j)E_j \| $. Considering orthogonal central projections,
one gets as a central carrier $Q_j$ of $E_j$ in $\Psi '$ as a
projection. Two $T$-algebras $\Psi $ and $\Theta $ are called
$T$-isomorphic, if an $\bf F$-linear multiplicative bijective
surjective mapping $\theta : \Psi \to \Theta $ exists continuous
together with its inverse mapping, \par $\theta (\alpha B)=\alpha
\theta (B)$, $ ~ \theta (AB) =\theta (A) \theta (B)$  and $\theta
(A^t)=[\theta (A)]^t$ for each $\alpha \in \bf F$, $A, B\in \Psi $.
Since $\theta $ and $\theta ^{-1}$ are continuous and
multiplicative, then $ \| \theta (A) \| = \| A \| $ is an isometry,
since \par $ \| \theta (\alpha A^n) \| = | \alpha |^n \| [\theta
(A)]^n \| \le |\alpha |^n \| \theta (A) \| ^n $ \\ for each $A\in
\Psi $ and $\alpha \in \bf F$. The $T$ algebras $ \bar{\Psi }E_j$
and $\bar{\Psi }Q_j$ are $T$-isomorphic, since $E^j, ~ I-E_j, ~ Q_j,
~ I-Q_j$ and $E_j^t, ~ (I-E_j)^t= I-E_j^t, ~ Q_j, ~ I-Q_j^t\in \Psi
$, where $\Psi $ is topologically complete. Then $ \| (\bar{D} A_j)
E_j \| = \| (\bar{D} A_j) Q_j \| = \| \bar{D} (A_j Q_j) \| \le \|
\bar{D} \| \| A_jQ_j \| = \| \bar{D} \| \| A_jE_j \| .$ From Theorem
10 it is known that $\bar{D}$ annihilates the center $Z(\bar{\Psi
})$ of $\bar{\Psi }$. Therefore, \par $ \| D_1 ( A_1E_1+...+A_nE_n)
\| \le \| \bar{D} \|
\max_j \| A_jE_j \| = \| \bar{D} \| \| A_1E_1+...+A_nE_n \| $, \\
consequently, $D_1$ is bounded. Thus $D_1$ has a bounded extension
being a derivation $D_1: \Upsilon \to \bar{\Upsilon }$. In view of
Lemma 12 it has a continuous extension $\bar{D}_1$ defined on
$\bar{\Upsilon }$.
\par On the other hand, $\bar{\Upsilon }$ is a $T$-algebra containing
$\bar{\Psi }$ and $\Xi $, since it contains $\cal P$, the projection
lattice of $\Xi $, hence $\Upsilon '\subset \Psi '$ and $\Upsilon '$
commutes with $\Xi $. But $\Xi $ is a maximal commutative subalgebra
in $\Psi '$, we get $\Upsilon ' = \Xi $.
\par Recall that a vector $x\in X$ is topologically cyclic relative to the action of $\Psi $
for a closed linear subspace $Y$ over $\bf F$ if $\bar{\Psi }x = \{
Ax: ~ A\in \bar{\Psi } \} $ is everywhere dense in $Y$. A subspace
$Y$ is called invariant relative to $\bar{\Psi }$, if $AY\subset Y$
for each $A\in \bar{\Psi }$. A closed linear subspace $Y$ in $X$
over $\bf F$ is called topologically irreducible relative to
$\bar{\Psi }$, if $Y$ is invariant relative to $\bar{\Psi }$ and
each non zero vector $x\in Y\setminus \{ 0 \} $ is topologically
cyclic relative to $\bar{\Psi }$. If $Y$ is a topologically
irreducible subspace, it has and orthocomplement $X\ominus Y$. So
$X\ominus Y$ has another topologically invariant subspaces and the
process can be done by transfinite induction (see \cite{eng}).
Therefore, the sum of all topologically irreducible subspaces in $X$
relative to $\bar{\Psi }$ is everywhere dense in $X$.
\par For any topologically irreducible subspace $Y$ relative to $\bar{\Psi }$
consider the restriction $\bar{\Psi }|_Y = \{ A|_Y: ~ A\in \bar
{\Psi } \} $. Since $\bar{D}_1A\in \bar{\Psi }$ for each $A\in \bar
{\Psi } $, the subspace $Y$ is invariant relative to
$\bar{D}_1\bar{\Psi }$ also. The algebra $\Psi $ and the Banach
space $X$ are over the spherically complete field $\bf F$. Take an
(ortho)projection $P$ from $X$ onto a finite dimensional over $\bf
F$ subspace $PX$ of a topologically irreducible subspace $Y$. This
induces the finite dimensional over $\bf F$ subalgebra $P\Psi P = \{
PAP: ~ A\in \Psi \} $. Then the differentiations $PDP: P\Psi P\to
P\Psi P$ and $PD_1P: P\Psi P\to P\Psi P$ act on it. \par Let $J_P$
be the center of $P\Psi P$. Then the differentiation operator $PDP$
annihilates $J_P$ due to Theorem 10 and hence $PD_1P$ annihilates
$J_P\cap \Upsilon $, consequently, $PDP$ and $PD_1P$ are defined on
the quotient algebras $(P\Psi P)/J_P$ and $(P\Upsilon P)/J_P$
correspondingly. Introduce on $(P\Psi P)/J_P$ the Lie algebra $\Psi
_P$ structure by $[A,B]=AB-BA$ for each $A, B\in (P\Psi P)/J_P$.
Traditionally $ad ~ B$ denotes $ad ~ B(A) = [B,A]$ for each $A\in
L(X)$. The latter Lie algebra $\Psi _P$ is non degenerate, i.e. has
a non degenerate Killing form $tr (ad A ~ ad B)$, where $(ad A)(E)=
[A,E]$ for each $A, E\in \Psi _P$. Then $PDP$ is the differentiation
of the Lie algebra $\Psi _P$ so that $PDP[A,B]= [PDPA,B]+[A,PDPB]$
and analogously $PD_1P$ is the differentiation of $\Upsilon _P$. In
view of Theorem 1.5.8 \cite{gotogrosb} the Lie algebra $\Psi _P$ is
complete, i.e. its center is zero and each its differentiation is
internal, $der (\Psi _P)=ad (\Psi _P)$, also $\Upsilon _P$ is
complete. Thus $PDP$ and $PD_1P$ are internal derivations of $P\Psi
P$ and $P\Upsilon P$ respectively.
\par Particularly, if $\bf
F$ is a locally compact field take $G_{\alpha }={\bf F}$. Generally
we consider a family $ \{ {\bf G}_{\alpha }: \alpha \in \mu \} $ of
locally compact subfields such that $\overline{\bigcup_{\alpha \in
\mu } {\bf G}_{\alpha }} = {\bf F}$. Since $\bf F$ is spherically
complete and ${\bf G}_{\alpha }$ is locally compact, then ${\bf
G}_{\alpha }$ is spherically complete. This family of subfields is
naturally directed by inclusion which induces a direction on $\mu $
such that $\alpha \le \beta $ if and only if ${\bf G}_{\alpha
}\subset {\bf G}_{\beta }$. Consider $\Psi $ over ${\bf G}_{\alpha
}$ and denote it by $\Psi _{\alpha }$. In view of Alaoglu-Bourbaki's
theorem (see \S 9.202 \cite{nari}) each bounded closed ball $B((\Psi
_{\alpha })',z,r)$ of radius $0<r<\infty $ and containing $z$ in
$(\Psi _{\alpha })'$ is weak-operator closed, since ${\bf G}_{\alpha
}$ is a locally compact field. \par From the proof above it follows
that $der (\Psi _{\alpha ,P})=ad (\Psi _{\alpha ,P})$ for each
$\alpha \in \mu $ and $P$ as above on $X_{\alpha }$, where
$X_{\alpha }$ is the Banach space $X$ considered over ${\bf
G}_{\alpha }$. The set ${\cal P}_{\alpha }$ of projections $P$ on
$X_{\alpha }$ is also directed by $P\le Q$ if and only if
$P(X_{\alpha })\subset Q(X_{\alpha })$. There are natural connecting
continuous $G_{\alpha }$-linear mappings $\pi ^{\beta }_{\alpha }:
X_{\beta } \to X_{\alpha }$ for each $\alpha \le \beta \in \mu $.
Put $B_{\alpha }$ to be the projective limit $B_{\alpha
}=\frac{\lim}{\leftarrow {\cal P}_{\alpha } } B_{\alpha ,P} $ which
exists in $(\Psi _{\alpha })'$. Then we put
$B=\frac{\lim}{\leftarrow \mu } B_{\alpha }$. These projective
limits exist relative to the weak-operator topology due to
Proposition 2.5.6 and Corollary 2.5.7 \cite{eng}. This operator $B$
is $\bf F$ linear, since it is ${\bf G}_{\alpha }$ linear on
$X_{\alpha }$ for each $\alpha $ and $\overline{\bigcup_{\alpha \in
\mu } {\bf G}_{\alpha }} = {\bf F}$.
\par Considering all possible topologically invariant subspaces and
all (ortho)projections $P$ with finite dimensional over $\bf F$
ranges one gets due to Theorem 15, that $\bar{D_1}$ is the internal
derivation of $\bar{\Upsilon }$, since the family of all finite
dimensional over $\bf F$ subalgebras $P\Upsilon P$ is everywhere
dense in $\bar{\Upsilon }$ relative to the weak-operator topology.
Then $D=ad B$ on $\Psi $ for some $B\in \Xi '$, since
$BT-TB=\bar{D}_1(T)=D(I)T=0$ for each $T\in \cal P$.

\par {\bf 18. Definition.} A derivation $D$ of an algebra $\Psi $ is
called inner, if $D=ad B|_{\Psi }$ for some element $B\in \Psi $ of
this algebra.

\par {\bf 19. Lemma.} {\it Each derivation $ad B$ of a $T$ algebra
$\Psi $ induces a derivation of $\Psi '$. A derivation $ad B$ of
$\bar{\Psi }$ is inner if and only if it induces an inner derivation
of $\Psi '$.}
\par {\bf Proof.} For every $A\in \Psi $ and $T\in \Psi '$ one gets
$(BT-TB)A-A(BT-TB)=BTA-TBA-ABT+ATB= (BA-AB)T - T(BA-AB)=0$, since
$[B,A]\in \Psi $. In the case when $ad B$ induces an inner
derivation of $\bar{\Psi }$ so that $ad B= ad E$ on $\bar{\Psi }$
with $E\in \bar {\Psi }$ this implies that $(B-E)$ commutes with
$\bar{\Psi }$. Therefore, $(B-E)\in \Psi '$. The inclusion $E\in
\bar{\Psi }$ implies that $ad (B-E)=ad B$ on $\Psi '$. That is $ad
B$ induces an inner derivation of $ \Psi '$.

\par {\bf 20. Definitions.} Suppose that $X$ is a Banach space over
a field $\bf F$ and $P$ is a projection on $X$, $~ P: X\to X$, and
$\Psi $ is a $W^t$ subalgebra in $L(X)$, $~P\in \Psi $. A projection
$P$ is called cyclic in $\Psi $ (or under $\Psi '$), if $PX=cl_X ~
span_{\bf F}\Psi 'x$ for some vector $x\in X$, where $span_{\bf F} U
:= \{ y\in X: ~ y = b_1x_1+...+b_nx_n; ~ b_1,...,b_n\in {\bf F},
~x_1,...,x_n\in X \} $, $cl_X U$ denotes the closure of a subset $U$
in $X$ relative to the norm topology. Such vector $x$ is called a
generating vector under $\Psi '$.
\par An orthoprojection $P$ in $\Psi $ over a spherically complete
field $\bf F$ is called countably decomposable relative to $\Psi $,
if every orthogonal family of non zero suborthoprojections of $P$ in
$\Psi $ is countable. When the unit operator $I$ is countably
decomposable relative to $\Psi $, one says that the $W^t$ algebra
$\Psi $ is countably decomposable.

\par {\bf 21. Lemma.} {\it Let $P$ be a central (ortho)projection in a $W^t$ algebra
$\Psi $ over a spherically complete field $\bf F$. This projection
$P$ is the central carrier of a cyclic projection in $\Psi $ if and
only if $P$ is countably decomposable relative to the center $Z(\Psi
)$ of $\Psi $. Moreover, a cyclic projection in $\Psi $ is countably
decomposable; two projections $P$ and $Q$ with the same generating
vector in $\Psi $ and $\Psi '$ have the same central carrier.}
\par {\bf Proof.} Consider a central projection $T$ in $\Psi $ with
generating vector $x\in X$ and $P=C_T$. Consider the case when there
are orthogonal families $P_{\alpha }$ and $T_{\beta }$ of
(ortho)projections in $Z(\Psi )$ and $\Psi $ respectively contained
in $P$ and $T$ correspondingly. The field $\bf F$ is spherically
complete and the Banach space $X$ is isomorphic with $c_0(\omega
,{\bf F})$ for some set $\omega $. Each closed linear subspace in
$X$ has an orthonormal basis which can be completed to an
orthonormal basis in $X$. If $y\in X$, then there are convergent
series $y=\sum_{\alpha }P_{\alpha }y$ and $y= \sum_{\beta } T_{\beta
}y$, where $P_{\alpha } y\perp P_{\beta }y$ and $T_{\alpha } y\perp
T_{\beta }y$ are orthogonal in the non-archimedean sense for each
$\alpha \ne \beta $. The convergence of these series is equivalent
to that for each $\epsilon
>0$ sets $\{ \alpha : ~ \| P_{\alpha }y \| > \epsilon \} $ and $ \{
\beta : ~ \| T_{\beta }y \| > \epsilon \} $ are finite. Thus these
series may have only countable sets of non zero additives. \par When
$T_{\beta }y=0$ the equalities $ \{ 0 \} = cl_X~ span_{\bf F} \Psi
'T_{\beta }y = cl_X span_{\bf F} T_{\beta } \Psi 'y$ and $T_{\beta
}Ty=T_{\beta }y$ are valid, if $P_{\alpha }y=0$ analogously
$P_{\alpha }y=0$. That is $P_{\alpha }Py=P_{\alpha }y=0$ due to the
equivalence of conditions $(i)$ and $(ii)$ in Section 17. Thus the
families $\{ P_{\alpha } \} $ and $ \{ T_{\beta } \} $ have at most
countable subsets of non zero elements, consequently, $P$ and $T$
are countably decomposable.
\par On the other hand, if $P$ is countably decomposable and $ \{
P_n \} $ is a countable set of projections cyclic under $(Z(\Psi
))'$ with generating vectors $x_n$ of unit norm and with sum $\vee_n
P_n=P$. The field $\bf F$ is of zero characteristic and contains the
$p$-adic field ${\bf Q}_p$ for some prime number $p$. Take the
vector $x=\sum_n p^n x_n$, where $n\in {\bf N}$. This sum or series
converges, since $ \| p^n x_n \| = \| x_n \| p^{-n}$ up to an
equivalence of norms on $\bf F$. Therefore the equality is valid $
cl_X ~ span_{\bf F} (Z(\Psi ))'x =PX$, since $cl_X ~ span_{\bf F}
(Z(\Psi ))'x$ contains $cl_X ~ span_{\bf F} (Z(\Psi ))'P_nx= cl_X ~
span_{\bf F} (Z(\Psi ))'x_n=P_nX$ for each $n$. Putting $T$ to be an
projection from $X$ onto $cl_X ~ span_{\bf F} (\Psi )'x$ one gets
$T\subseteq P$, since $\Psi '\subseteq (Z(\Psi ))'$, that is $PT=T$.
Suppose that $Q\in Z(\Psi )$ and $QT=T$. This implies that $Qx=x$
and $cl_X ~ span_{\bf F}(Z(\Psi ))' x = cl_X ~ span_{\bf F}(Z(\Psi
))' Qx = cl_X ~ span_{\bf F}Q(Z(\Psi ))' x$, consequently, $P=QP$.
This means that $P=C_T$ with $T$ cyclic in $\Psi $. Therefore, the
projection $P$ is the central carrier of $cl_X ~ span_{\bf F} \Psi
x$.

\par {\bf 22. Theorem.} {\it If $\Psi $ is a $W^t$ algebra on a Banach space
over a spherically complete field $\bf F$ and $D$ is a derivation of
$\Psi $, then $D$ is inner.}
\par {\bf Proof.} In view of Theorem 17 a derivation $D$ has the
form $D=ad B|_{\Psi }$ for some bounded linear operator $B\in L(X)$.
Then $-(BA^t-A^tB)^t=B^tA-AB^t\in \Psi $ for each $A\in \Psi $,
consequently, the mapping $ad B^t: \Psi \to \Psi $ is also the
differentiation of $\Psi $. Therefore, $ad (B+B^t)$ and $ad (B-B^t)$
are derivations of $\Psi $. If each of these derivations $ad
(B+B^t)$ and $ad (B-B^t)$ is inner, then $ad B$ is inner as well.
Mention that the operator $ad (\lambda I+B)$ is the derivation
together with $ad B$ for each $\lambda \in \bf F$. In accordance
with Theorem 10 $B\Phi '$ is the center of the $W^t$ algebra $\Psi
$, where $\Phi =Z(\Psi )$. \par If $\{ P_{\beta }: ~ \beta \in
\Lambda \} $ is a family of projections on $X$ so that its sum
$I=\sum_{\beta \in \Lambda } P_{\beta }$ is the unit operator and
$ad B|_{\Psi P_{\beta }} = ad E_{\beta }|_{\Psi P_{\beta }}$ for
every $\beta $ and $\sup_{\beta \in \Lambda } \| E_{\beta } \|
<\infty $, where $\Lambda $ is a suitable set, $E_{\beta } \in \Psi
P_{\beta }$, then $ad B|_{\Psi } = ad E|_{\Psi }$ for $E=\sum_{\beta
\in \Lambda } E_{\beta }$.
\par Take $Q_{\alpha }$ a cyclic projection under $(Z(\Psi ))'$
for each $\alpha $. It is sufficient to prove this assertion for
countably decomposable center $Z(\Psi )$ due to Lemma 21. For this
one takes a cyclic projection $T$ in $\Psi '$ with central carrier
$I$ considering the faithful representation $\Psi T$ of $\Psi $ on
$T(X)$. The commutant is $T\Psi 'T$ and so it is sufficient to
consider that $ \Psi '$ is countably decomposable.
\par Let $\bf G$ be a locally compact field contained in
$\bf F$ and consider the spherically complete field $\bf F$ as the
Banach space over $\bf G$ isomorphic with $c_0(\omega ,{\bf F})$ for
some set $\omega $ (see \S 21). Then the Banach space $X$ over $\bf
F$ has the structure of the Banach space $X_{\bf G}$ over $\bf G$ as
well. To each operator bounded linear operator $A\in L(X)$ a bounded
operator $A_{\bf G}\in L(X_{\bf G})$ corresponds. Due to
Alaoglu-Bourbaki's theorem (see \S 9.202 \cite{nari}) a closed
bounded ball $B(X_{\bf G},x,r):= \{ y\in X_{\bf G}: ~ \| y-z \| \le
r \} $ in $X_{\bf G}$ is weakly compact and a bounded closed ball
$B(L(X_{\bf G}),A,r):= \{ C\in L(X_{\bf G}): ~ \| C-A \| \le r \} $
in $L(X_{\bf G})$ is compact relative to the weak operator topology,
where $0<r<\infty $. Therefore, $B(L(X_{\bf G}),A,r)\cap \Psi _{\bf
G}$ is also compact relative to the weak operator topology, where
$\Psi _{\bf G}$ is the $W^t$ algebra $\Psi $ considered over the
field $\bf G$, i.e. by narrowing the field from $\bf F$ to $\bf G$
so that $\Psi _{\bf G}\subset L(X_{\bf G})$.
\par A system of algebras $ \{ \Psi _P: ~ P \in {\cal P} \} $ and a
family of locally compact subfields $\{ {\bf G}_{\alpha }: \alpha
\in \mu \} $ from \S 17 gives rise to the projective limit
decomposition of each operator $A\in \Psi $ or $E\in \Psi '$ and for
the differentiation operator $D$ as well, since $\Psi = \bar{\Psi }$
by the conditions of this theorem. Finally, from Proposition 2.5.6
and Corollary 2.5.7 \cite{eng} the assertion follows.


\begin{thebibliography}{99}

\bibitem{diarlud02} B. Diarra, S.V. Ludkovsky. "Spectral integration and spectral theory
for non-archimedean Banach spaces"// Intern. J. of Mathem. and
Mathem. Scinces {\bf 31: 7 } (2002), 421-442.

\bibitem{eng} R. Engelking. "General topology" (Mir: Moscow, 1986).

\bibitem{gotogrosb} M. Goto, F.D. Grosshans. "Semisimple Lie
algebras" (Marcel Dekker, Inc.: New York, 1978).

\bibitem{kadannm66} R.V. Kadison. "Derivations of operator
algebras"// Annals of Mathem. {\bf 83: 1} (1966), 280-293.

\bibitem{kadannmat57} R.V. Kadison. "Unitary invaraints for representations of operator
algebras"// Annals of Mathem. {\bf 66: 2} (1957), 304-379.

\bibitem{kadringb} R.V. Kadison, J.R. Ringrose. "Fundamentals of the theory of operator algebras"
(Academic Press: New York, 1983).

\bibitem{kakolinm94} J. Kakol. "Remarks on spherical completeness of
non-archimedean valued fields"// Indag. Mathem. {\bf 5: 3 } (1994),
321-323.

\bibitem{losannm2008} V. Losert. "The derivation problem for
group algebras"// Annals of Mathem. {\bf 168: 1} (2008), 221-246.

\bibitem{ludjmsqim05} S.V. Ludkovsky. "Quasi-invariant and
pseudo-differentiable measures with values in non-archimedean fields
on a non-archimedean Banach space"// J. Mathem. Sci. {\bf 128: 6}
(2005), 3428-3460.

\bibitem{nari} L. Narici, E. Beckenstein. "Topological vector
spaces"  (Marcel Dekker, Inc.: New York, 1985).

\bibitem{putmatann72} M. van der Put. "Difference equations over
$p$-adic fields"// Math. Ann. {\bf 198} (1972), 189-203.

\bibitem{roo} A.C.M. van Rooij. "Non-Archimedean functional
analysis" (Marcel Dekker, Inc.: New York, 1978).


\bibitem{sakaiannm66} S. Sakai. "Derivations of $W^*$-algebras"//
Annals of Mathem. {\bf 83: 1} (1966), 273-279.


\bibitem{schikb} W.H. Schikhof. "Ultrametric calculus"
(Cambridge Univ. Press: Cambridge, 1984).

\bibitem{wei} A. Weil. "Basic number theory" (Springer: Berlin,
1973).

\end{thebibliography}
\end{document}